\documentclass[portuges,12pt,letter]{article}
\usepackage[centertags]{amsmath}
\usepackage{amsfonts}
\usepackage{newlfont}
\usepackage{amscd}
\usepackage{graphics}
\usepackage{epsfig}
\usepackage{indentfirst}
\usepackage{amsxtra}
\usepackage[latin1]{inputenc}
\usepackage{amssymb, amsmath}
\usepackage{amsthm}
\usepackage[mathscr]{eucal}

\newtheorem{thm}{Theorem}[section]

\newtheorem{remark}[thm]{Remark}

\author{Fabio Silva Botelho \\ Department of Mathematics \\  Federal University of Santa Catarina, UFSC \\
Florian\'{o}polis, SC - Brazil}

\title{\bf  A numerical method for an inverse optimization problem through the generalized method of lines}

\begin{document}
\maketitle

\abstract{ This article develops a solution for an inverse problem through the generalized method of lines. We consider a Laplace equation on a domain with internal and external boundaries with standard Dirichlet boundary conditions. Also, we specify a third non-homogeneous Newmann type boundary condition for the external boundary, and consider the  problem of finding the optimal shape for the internal boundary such that all the prescribed boundary conditions are satisfied. The novelty here presented is the application of the generalized method of lines as a tool to compute a solution for such an inverse optimization  problem.}

\section{Introduction} In this article we develop a numerical method to compute the solution of an inverse problem through the generalized method of lines.

More specifically, we consider a Laplace equation on a domain $\Omega \subset \mathbb{R}^2$ with an internal boundary denoted by $\partial \Omega_0$ and
an external one denoted by $\partial \Omega_1.$ We prescribe boundary conditions for both $\partial \Omega_0$ and $\partial \Omega_1$ and a third boundary condition for $\partial \Omega_1$ (the external one) and consider the problem of finding the optimal shape for the internal boundary $\partial \Omega_0$ for which such a third boundary condition for the external boundary is satisfied.

The idea is to discretize the domain in lines (in fact curves) and write the solution of Laplace equation on these lines as functions of the unknown internal
boundary shape, through the generalized method of lines.

The second step is to minimize a functional which corresponds to the $L^2$ norms of Laplace equation and concerning third boundary condition.
\begin{remark} About the references, this and many other similar problems are addressed in \cite{770}. The generalized method of lines has been originally
introduced in \cite{901}, with additional results in \cite{909,12a}. Moreover, about finite differences schemes we would cite \cite{103}. Finally, details on the function spaces here addressed may be found in \cite{1}.
\end{remark}
\section{ The mathematical  description of the main problem}

At this point we start to describe mathematically our main problem.

Let $\Omega \subset \mathbb{R}^2$ be a bounded, closed and  connected  set defined by

$$\Omega=\{(r,\theta) \in \mathbb{R}^2\;:\; r(\theta) \leq r \leq R\;:\; 0 \leq \theta \leq 2\pi\}.$$

Consider a Laplace equation and concerning boundary conditions expressed by

\begin{equation} \left\{
\begin{array}{ll}
 \nabla^2 u=0,& \text{ in } \Omega,
 \\
 u=u_o, & \text{ on } \partial \Omega_0,
 \\
 u=u_f, & \text{ on } \partial \Omega_1,
 \\
 \nabla u \cdot \mathbf{n}=w, &\text{ on } \partial \Omega_1,\end{array} \right.\end{equation}
 where $u_o, u_f \text{ and } w \in C^2([0,2\pi])$ are known periodic functions with period $2\pi,$
  $\mathbf{n}$ denotes the outward normal field to $\partial \Omega,$
 $$\partial \Omega_0=\{(r(\theta),\theta) \;:\; 0 \leq \theta \leq 2\pi\},$$
 $$\partial \Omega_1=\{(R, \theta)\;:\; 0 \leq 0 \leq 2\pi\},$$
 and $R>0$.

The main idea is to discretize the domain in lines and through the generalized method of lines  to obtain the solution $u_n(r(\theta))$ on each line $n$ as a function
of $r(\theta).$ The final step is to compute the optimal $r(\theta)$ in order to minimize the cost functional $J(r(\theta))$ defined by
$$J(r(\theta))=\| \nabla^2 u\|_{2,\Omega}^2+K \|\nabla u \cdot \mathbf{n} -w\|_{2,\partial \Omega_1}^2,$$ where $K>0$ is an appropriate constant
to be specified.

\section{About the generalized method of lines and the  main result}

At this point we start to describe the main result.

Consider firstly a Laplace equation in polar coordinates, that is,

$$\frac{\partial^2 u}{\partial r^2}+\frac{1}{r} \frac{\partial u}{\partial r}+\frac{1}{r^2}\frac{\partial^2 u}{\partial \theta^2}=0.$$

In order to apply the generalized method of lines, we a define a variable $t$, through the equation
$$ t=\frac{r-r(\theta)}{R-r(\theta)}$$  so that $t \in [0,1]$ in the set $\Omega$ previously specified.

Hence, denoting $u(r,\theta)=\overline{u}(t,\theta)$,  we obtain,
\begin{eqnarray}
\frac{\partial u}{\partial \theta}&=& \frac{\partial \overline{u}}{\partial t} \frac{\partial t}{\partial \theta} +\frac{\partial \overline{u}}{\partial \theta}
\nonumber \\ &=& \frac{\partial \overline{u}}{\partial t}\left(\frac{t-1}{R-r(\theta)}\right) r'(\theta)+\frac{\partial \overline{u}}{\partial \theta} \nonumber \\ &=&
\frac{\partial \overline{u}}{\partial t}(f_1(\theta)+t f_2(\theta))+\frac{\partial \overline{u}}{\partial \theta},\end{eqnarray}
where
$$f_1(\theta)=\frac{-r'(\theta)}{R-r(\theta)},$$
and
$$f_2(\theta)=\frac{r'(\theta)}{R-r(\theta)}.$$

Also,
\begin{eqnarray}
\frac{\partial u}{\partial r}&=&\frac{\partial \overline{u}}{\partial t}\frac{\partial t}{\partial r}
\nonumber \\ &=& \frac{\partial \overline{u}}{\partial t}\frac{1}{R-r(\theta)},\end{eqnarray}
so that
\begin{eqnarray}\frac{\partial^2 u}{\partial r^2}&=&\frac{\partial^2 \overline{u}}{\partial t^2}\left(\frac{\partial t}{\partial r}\right)^2
\nonumber \\ &=& \frac{\partial^2 \overline{u}}{\partial t^2}\frac{1}{(R-r(\theta))^2} \nonumber \\ &=& f_3(\theta)\frac{\partial^2 \overline{u}}{\partial t^2},\end{eqnarray}
where $$f_3(\theta)=\frac{1}{(R-r(\theta))^2}.$$

Moreover, we may also obtain
\begin{eqnarray}
\frac{\partial^2 u}{\partial \theta^2}&=& f_4(t,\theta)\frac{\partial^2 \overline{u}}{\partial t^2}
+f_5(t,\theta)\frac{\partial \overline{u}}{\partial t}\nonumber \\ &&+f_6(t,\theta)\frac{\partial^2 \overline{u}}{\partial t \partial \theta}
+ \frac{\partial^2 \overline{u}}{\partial \theta^2},
\end{eqnarray}
where,
$$f_4(t,\theta)=(f_1(\theta)+tf_2(\theta))^2,$$
$$f_5(t,\theta)=f_1'(\theta)+t f_2'(\theta)+f_2(\theta)(f_1(\theta)+tf_2(\theta)),$$
$$f_6(t,\theta)=2(f_1(\theta)+t f_2(\theta)).$$

Thus, for the new variables $(t,\theta)$, dropping the bar in $\overline{u}$, the Laplace equation is equivalent to
$$\frac{\partial^2u}{\partial t^2}+f_7(t,\theta) \frac{\partial u}{\partial t}+f_8(t, \theta) \frac{\partial^2u}{\partial t \partial \theta}
+f_9(t,\theta)\frac{\partial^2 u}{\partial \theta^2}=0,$$  in $\hat{\Omega}$ where
$$\hat{\Omega}=\{(t,\theta) \in \mathbb{R}^2 \;:\; 0\leq t \leq 1,\; 0\leq \theta \leq 2\pi\},$$
$$r=t(R-r(\theta))+r(\theta),$$
$$f_0(t,\theta)=\frac{1}{(R-r(\theta))^2}+\frac{f_4(t,\theta)}{r^2},$$
$$f_7(t,\theta)= \frac{\tilde{f}_7(t,\theta)}{f_0(t,\theta)},$$
$$f_8(t,\theta)=\frac{\tilde{f}_8(t,\theta)}{f_0(t,\theta)},$$
$$f_9(t,\theta)=\frac{\tilde{f}_9(t,\theta)}{f_0(t,\theta)}$$
and where
$$\tilde{f}_7(t,\theta)=\left(\frac{1}{r(R-r(\theta))}+\frac{f_5(t,\theta)}{r^2}\right)$$
$$\tilde{f}_8(t,\theta)=\frac{f_6(t,\theta)}{r^2 }$$
and
$$\tilde{f}_9(t,\theta)=\frac{1}{r^2}.$$

So, discretizing the interval $[0,1]$ into $N \in \mathbb{N}$ pieces, that is defining, $$t_n=\frac{n}{N},\; \forall n \in\{0,\ldots, N\}$$
and $d=1/N$, in partial finite differences, the concerning equation stands for
\begin{eqnarray}\label{eq1}
&&\frac{u_{n+1}-2u_n+u_{n-1}}{d^2}+f_7(t_n,\theta) \left(\frac{u_n-u_{n-1}}{d} \right)
\nonumber \\ &&+ f_8(t_n,\theta)\frac{\partial}{\partial \theta}\left(\frac{u_n-u_{n-1}}{d} \right)+
f_9(t_n,\theta) \frac{\partial^2 u_n}{\partial \theta^2}=0,\end{eqnarray}
$\forall n \in \{1,\ldots,N-1\},$ where
$$u_0=u_o(\theta)$$ and $$u_N=u_f(\theta).$$

At this point we describe how to obtain the general expression for $u_n$ corresponding to the line $n$, through the generalized method of lines.

For $n=1$ in (\ref{eq1}), we have
\begin{eqnarray}
&&u_2-2u_1+u_0+f_7(t_1,\theta) (u_1-u_0) d
\nonumber \\ &&+ f_8(t_1,\theta)\frac{\partial (u_1-u_0)}{\partial \theta}d+
f_9(t_1,\theta) \frac{\partial^2 u_1}{\partial \theta^2}d^2=0,\end{eqnarray}
so that
$$u_1=T_1(u_2,u_1,u_0),$$
where \begin{eqnarray}
T_1(u_2,u_1,u_0)&=&\left(u_2+u_1+u_0+f_7(t_1,\theta) (u_1-u_0) d \right.
\nonumber \\ &&\left.+ f_8(t_1,\theta)\frac{\partial (u_1-u_0)}{\partial \theta}d+
f_9(t_1,\theta) \frac{\partial^2 u_1}{\partial \theta^2}d^2 \right)/3.\end{eqnarray}

To solve this equation we apply the Banach fixed point theorem as it follows:

\begin{enumerate}
\item Set $u_1^1=u_2.$
\item Define $$u_1^{k+1}=T_1(u_2,u_1^k,u_0),\; \forall k \in \mathbb{N}.$$
\item Obtain $$u_1=\lim_{k \rightarrow \infty} u_1^k\equiv F_1(u_2,u_0).$$
\end{enumerate}

Reasoning inductively, having
$$u_{n-1}=F_{n-1}(u_n,u_0)$$
also  from (\ref{eq1}), for the line $n$, we have

$$u_n=T_n(u_{n+1},u_n,u_0),$$
where \begin{eqnarray}
T_n(u_{n+1},u_n,u_0)&=&\left(u_{n+1}+u_n+F_{n-1}(u_n,u_0)+f_7(t_n,\theta) (u_n-F_{n-1}(u_n,u_0)) d \right.
\nonumber \\ &&\left.+ f_8(t_n,\theta)\frac{\partial (u_n-F_{n-1}(u_n,u_0))}{\partial \theta}d+
f_9(t_n,\theta) \frac{\partial^2 u_n}{\partial \theta^2}d^2 \right)/3.\end{eqnarray}

Again, to solve this equation, we apply the Banach fixed point theorem, as it follows:

\begin{enumerate}
\item Set $u_n^1=u_{n+1}.$
\item Define $$u_n^{k+1}=T_n(u_{n+1},u_n^k,u_0),\; \forall k \in \mathbb{N}.$$
\item Obtain $$u_n=\lim_{k \rightarrow \infty} u_n^k\equiv F_n(u_{n+1},u_0),\; \forall n \{2,\ldots,N-1\}.$$
\end{enumerate}

In particular, for $n=N-1$, we get
$$u_{N-1}=F_{N-1}(u_N,u_0)=F_{N-1}(u_f,u_0).$$

Similarly, for $n=N-2$, we obtain
$$u_{N-2}=F_{N-2}(u_{N-1},u_0),$$
and so on, up to finding
$$u_1=F_1(u_2,u_0).$$

This means that we have obtained a general expression $$u_n=F_n(u_{n+1},u_0)\equiv H_n(u_f,u_0,r(\theta)),\; \forall n \in \{1,\ldots,N-1\},$$

Anyway, we remark to properly run the software we have to make the approximation $t_n=t, \forall n  \in \{1, \ldots, N-1\}$. Thus, we have used the method
above described just to find a general expression for $u_n$, which is approximately given by (here $x$ stands for $\theta$)
\begin{eqnarray}u_n(x)&\approx& a_n[1] u_f(x)+a_n[2]u_0(x)+ a_n[3]\;f_7[t_n,x] u_f(x)+a_n[4]\;f_7[t_n,x] u_0(x) \nonumber \\ &&
+a_n[5]\;f_8[t_n,x] u_f'(x)+a_n[6]\;f_8[t_n,x] u_0'(x) \nonumber \\ &&+a_n[9]\;f_9[t_n,x] u_f''(x)+a_n[10]\;f_9[t_n,x] u_0''(x).
\end{eqnarray}
Indeed this a first approximation for the series representing the solution on each line obtained by considering terms up to order $d^2$ (in d).

Moreover the coefficients $a_n[k]$ and $r(\theta)$ expressed in finite differences are calculated through the minimization of
$J(r(\theta),\{a_n[k]\})$ given by
$$J(r(\theta),\{a_n[k]\})=\| \nabla^2 u\|_{2,\Omega}^2+K \|\nabla u \cdot \mathbf{n} -w\|_{2,\partial \Omega_1}^2,$$

\subsection{The numerical results}

We present numerical results for $R=30$, $K=250$;
$$u_o(x)=0.5\cos(\pi x)+0.8$$
$$u_f(x)=0.5 \sin(2\pi x)+1.0$$
and $$w(x)=0.3*(\cos(2\pi x)/2+1.0)$$

\begin{figure}
\centering \includegraphics [width=3in]{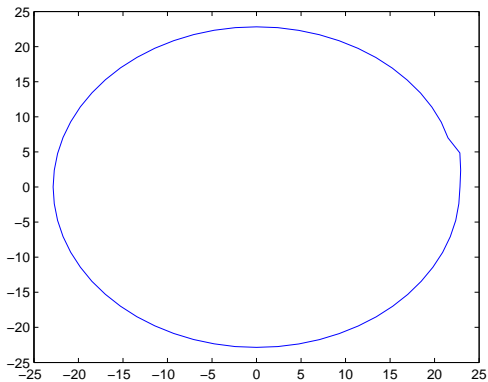}
\\ \caption{\small{ Optimal shape  for the internal boundary $\partial \Omega_0$ of $\Omega$.}}\label{shape-inverse-1}
\end{figure}

Please see figure \ref{shape-inverse-1}, for the optimal shape which minimizes $J$.

Finally, we have obtained,
$$\|(\nabla u)\cdot \mathbf{n}-w\|_\infty \approx 0.0167.$$
and
$$\|\nabla^2 u\|_\infty \approx 0.0111$$

These last two results indicate the method proposed has been successful to compute such a problem.

\section{Conclusion}
In this article we have used the generalized method of lines to obtain an approximate solution for an inverse optimization problem.

We emphasize the numerical results obtained indicate the numerical procedure proposed is an interesting possibility to compute
such a type of problem.

We also highlight the method here developed may be applied to a large class of similar models.

\end{document}